\documentclass[12pt,a4paper]{article}
\usepackage{amsmath,amsthm,amsfonts,amssymb,alltt}
\usepackage{graphics}
\usepackage{graphicx}
\usepackage{epsfig}
\usepackage{color}
\usepackage{longtable}

\theoremstyle{definition}

\theoremstyle{remark}

\numberwithin{equation}{section}

\begin{document}

\author{Semyon Yakubovich$^1$, Piotr Drygas$^2$, Vladimir Mityushev$^3$}

\title{Closed-form evaluation of 2D static lattice sums}
\maketitle
{
\noindent$^1$Department of Mathematics,  Faculty of Sciences,  University of Porto,  Campo Alegre st., 687,  4169-007 Porto,    Portugal, \\syakubov@fc.up.pt
\\
$^2$Department of Differential Equations and Statistics, 
Faculty of Mathematics and Natural Sciences, 
University of Rzeszow,
 Pigonia 1, 
35-959 Rzeszow, Poland, \\drygaspi@ur.edu.pl
\\
$^3$Institute of Computer Sciences and Computer Methods, Pedagogical University, ul.Podchorazych 2, Krakow, 30-084, Poland,  \\mityu@up.krakow.pl
\\
$^3$The corresponding author}

%


\section{Introduction}
The mathematical questions of convergence, numerically effective algorithm and closed-form evaluation of the lattice sums were discussed in the fundamental book Borwein et al \cite{Borwein} and works cited therein. The present paper is devoted to closed-form evaluation of the conditionally convergent 2D lattice sums. One of them, $S_2$, defined by \eqref{eq:s5}, was considered by Lord Rayleigh  \cite{Rayleigh}. Numerically effective series for Rayleigh's sum  based on the elliptic functions are outlined in Borwein et al \cite[Sec.3.2]{Borwein}.  McPhedran et al  \cite{McPhedran}, {Movchan et al  \cite{Movchan} and Greengard et al  \cite{Greengard} developed the Rayleigh method to elastostatic  (see lattice sum \eqref{eq:g100}) and elastodynamic problems having paid the main attention to computationally convenient and accurate expressions for the lattice sums constructed for the square array. 

The effective properties of unidirectional fibrous  composites  can be expressed in terms of the series in concentration $f$. The famous Clausius-Mossotti approximation also known as the Maxwell formula \cite[Ch. 10]{Milton} is valid in the first order approximation. The second order approximation includes the value of $S_2$. In order to shortly describe these approximations following Rayleigh we consider a doubly periodic rectangular array of disks of conductivity $\lambda_1$ embedded in matrix of conductivity $\lambda$. Let $\lambda_{xx}$ and $\lambda_{yy}$ be the principal components of the effective conductivity tensor.  Then \cite{Rayleigh}, \cite{MR2001}, \cite{MR2013}    
\begin{align}
\label{eq:CMA}
\frac{\lambda_{xx}}{\lambda} &= 1+2\rho f+2\rho^2 f^2 \frac{S_2}{\pi} + O((|\rho| f)^3),\\
\frac{\lambda_{yy}}{\lambda} &= 1+2\rho f+2\rho^2 f^2 \left(2-\frac{S_2}{\pi} \right)+ O((|\rho| f)^3),
\end{align}
where $\rho = \frac{\lambda_1-\lambda}{\lambda_1+\lambda}$ denotes the contrast parameter. For the square array, the medium becomes macroscopically isotropic, i.e., $\lambda_{e}=\lambda_{xx}=\lambda_{yy}$ and the above formulae becomes the Clausius-Mossotti approximation
\begin{equation}
\frac{\lambda_{e}}{\lambda} =1+2\rho f+2\rho^2 f^2  + O((|\rho| f)^3) = 
\frac{1+\rho f}{1-\rho f}+ O((|\rho| f)^3).
\label{eq:CMA2}
\end{equation}
Actually, the approximation in the right part of \eqref{eq:CMA2} holds up to $O((|\rho| f)^5)$ (see formula (28) from \cite{rylko2000} where the correction in the fifth order term should be $6S_4^2 \pi^{-2} \frac{f^5}{(1-f)^2}$).

The same rule holds for elastostatic problems \cite{arx}. An analytic formula for the macroscopic elastic constants must include lattice sums \eqref{eq:s5} and \eqref{eq:g100} in the second order term $O(f^2)$. Such formulae for elastic problems  are similar to \eqref{eq:CMA}-\eqref{eq:CMA2}. They are described in Sec.\ref{sec:5}. Therefore, analytical formulae for the lattice sums \eqref{eq:s5} and \eqref{eq:g100} have the fundamental applications in 2D composites. 

In the present paper, employing properties of the complete elliptic integrals of the first and second kind, we deduce closed-form formulae for the lattice sums \eqref{eq:s5},  \eqref{eq:g100} and other new formulae.  Applications to the effective properties of regular and random composites are discussed. 

\section{Eisenstein summation method and Rayleigh integral}
\label{sec:2}
Let $\mathbb Z$ and $\mathbb C$ denote the sets of integer and complex numbers, respectively, $i$ the imaginary unit. Consider a lattice \{$m \omega_1+n \omega_2 \in \mathbb C: m,n\; \mbox{run over}\; \mathbb Z\}$ determined by two fundamental translation vectors expressed by complex numbers $\omega_1$, $\omega_2$. Without loss of generality we assume that $\omega_1>0$ and  $\mathrm{Im}\tau>0$ where $\tau=\frac{\omega_2}{\omega_1}$. 
Let the area of the fundamental parallelogram be normalized to unity, hence, $\omega_1^2\mathrm{Im}\tau=1$. The main object of the present paper is the conditionally convergent lattice sums
\begin{equation}
S_2 =
\sideset{}{^e}\sum_{(m,n) \in \mathbb Z^2 \backslash \{(0,0)\}}
\frac{1}{(m \omega_1+n \omega_2)^2} = 
\mathrm{Im}\tau \sideset{}{^e}\sum_{(m,n) \in \mathbb Z^2 \backslash \{(0,0)\}} \frac{1}{(m+n \tau)^2}
\label{eq:s5}
\end{equation}
and
\begin{equation}
\label{eq:g100}
T_2 = 
\sideset{}{^e}\sum_{(m,n) \in \mathbb Z^2 \backslash \{(0,0)\}}
\frac{\overline{m \omega_1+n \omega_2}}{(m \omega_1+n \omega_2)^3} =
 \mathrm{Im}\tau \sideset{}{^e}\sum_{(m,n) \in \mathbb Z^2 \backslash \{(0,0)\}}
\frac{m+n \overline{\tau}}{(m+n \tau)^3},
\end{equation}
where the Eisenstein summation method \cite{Weil} is used 
\begin{equation}
\label{eq:eisen_sum}
\sideset{}{^e}\sum_{m,n}:=\lim_{M_2\rightarrow\infty}\sum_{n=-M_2}^{M_2}
\left(\lim_{M_1\rightarrow\infty}\sum_{m=-M_1}^{M_1}\right).
\end{equation}
Having used the summation method \eqref{eq:eisen_sum} Rayleigh deduced the formula
\begin{equation}
\label{eq:Ray}
S_2(\tau) = \pi^2\mbox{Im}\; \tau 
\left[\frac 13+2\sum _{m=1}^{\infty } \frac{1}{\sin^2(\pi  m \tau)} \right]
\end{equation}
for an rectangular array. It can be easily extended to other shapes of the fundamental cell. A physical justification of  the Eisenstein summation method  was presented in  \cite{McPhedran1978} and \cite{perrins}. A rigorous mathematical proof can be found in \cite{mit97}.
Rayleigh (1892) did not cite Eisenstein's result (1847) and addressed to Weierstrass' investigations (1856). Perhaps, it is related to that Eisenstein treated formally his series without uniform convergence introduced by Weierstrass. 

Moreover, Rayleigh \cite{Rayleigh} found the beautiful formula $S_2(i) = \pi$ where $\tau = i$ corresponds to the square array. The method of calculation was based on the reduction of the sum $S_2(i)$ to the integral
\begin{equation}
\label{eq:Ray1}
S_2(i)  = 2 \int_v^{\infty} dx \int_{-v}^v \frac{dy}{(x+iy)^2} =\pi.
\end{equation}
The integral over the central square $(-v,v) \times (-v,v)$ is eliminated in \eqref{eq:Ray1} since it vanishes. 
Further, the equality $S_2 = \pi$ was proved in \cite{perrins} for the hexagonal array (under the normalization of the area of the fundamental cell to unity).

The following formula was independently deduced in \cite{mitZAMM} 
\begin{equation}
\label{eq:mit1}
S_2(\tau)=\frac2{\omega_1} \zeta \left(\frac{\omega_1}2\right)
= \pi^2\mbox{Im}\; \tau 
\left[\frac 13-8\sum _{m=1}^{\infty } 
\frac{m \exp(2i\pi  m \tau)}{1-\exp(2i\pi  m \tau)} \right],
\end{equation}
where the $\zeta$-Weierstrass function is used. The equality $S_2(i) = \pi$ was proved by Legendre's identity. Though formulae \eqref{eq:Ray} and \eqref{eq:mit1} are similar a reduction of one to other can be justified only through the Eisenstein summation method applied to the $\zeta$-Weierstrass function.  

We now proceed to discuss the lattice sum \eqref{eq:g100} beginning form the relation \cite{Weil}
\begin{equation}
\label{eq:eps3}
\sum\limits_{m \in \mathbb Z \backslash \{0\}}\frac{1}{(m+\tau)^3}=\pi^3\frac{\cos (\pi\tau)}{\sin^3 (\pi \tau)} .
\end{equation}
Consider the general term of the series \eqref{eq:g100}
\begin{equation}
\label{eq:TT2}
\frac{m+n \overline{\tau}}{(m+n \tau)^3}=\frac{1}{(m+n \tau)^2} - 2 i \mbox{Im} \tau \frac{n}{(m+n \tau)^3}.
\end{equation} 
Substitution of \eqref{eq:eps3} and \eqref {eq:TT2}  into  \eqref{eq:g100} and use of \eqref{eq:s5} yields the computationally effective formula
\begin{equation}
\label{eq:numT}
T_2(\tau)= S_{2}(\tau)-4i\pi^3 (\mbox{Im}\;\tau)^2 \sum\limits_{n=1}^\infty n \frac{\cos (n\pi \tau)}{\sin^3 (n\pi \tau)}.
\end{equation} 

Surprisingly, that Rayleigh's integral gives a wrong result
\begin{equation}
\label{eq:Ray3}
T_2(i)  = 2 \int_v^{\infty} dx \int_{-v}^v \frac{x-iy}{(x+iy)^3} dy =2 
\end{equation}
though it formally corresponds to the Eisenstein summation method.  This is because the Rayleigh reduction to the integral is formal.  Moreover, by simple substitution we observe that iterated integrals \eqref{eq:Ray1}, \eqref{eq:Ray3} do not depend on $v$ and the corresponding double integrals diverge.  Therefore, only formula \eqref{eq:numT} was in our disposal to get numerical values of $T_2(\tau)$.    However, in the sequel we will give a rigorous proof of the closed-form formulae for the lattice sums  \eqref{eq:s5},   \eqref{eq:g100} on the imaginary axis and on the vertical lines ${\rm Re} \ \tau = \pm 1/2$, basing on the theory of the complete elliptic integrals of the first and second kind. 


\section{Closed-form formulae}
\label{sec:3}

Let $\mathbb{R},\ \mathbb{R}_+$ be the sets of real and real positive numbers, respectively.    Let $x \in \mathbb{R}_+$ be given by the formula
$$x\equiv x(k)= { K (k^\prime)\over K(k)},\  k \in (0, 1),\ k^\prime = \sqrt {1-k^2},\eqno(3.1)$$
where  $K(k)$ is the complete elliptic integral of the first kind \cite{apo}, \cite{erd}, Vol. II, \cite{wit} 

$$K(k)= \int_0^1 {dt\over \sqrt{ (1-t^2) (1- k^2 t^2)} }.\eqno(3.2)$$  
The parameter $k$ is called the elliptic modulus and $k^\prime$ is the complimentary modulus.   As we see the function $x$ as a function of  $k \in (0,1)$   is monotone decreasing and continuously differentiable bijective map $ x: (0,1) \to \mathbb{R}_+$.  Therefore any $x >0$ is uniquely defined by the corresponding modulus $k$.  The complete elliptic integral K(k) satisfies the Legendre relation

$$E(k)K (k^\prime) + E (k^\prime)K(k)-  K (k^\prime) K(k) = {\pi\over 2},\eqno(3.3)$$
where $E(k)$ is the complete elliptic integral of the second kind

$$E(k)= \int_0^1 \sqrt{{ 1- k^2 t^2\over 1-t^2} } dt.\eqno(3.4)$$  
Its derivative can be calculated by formula 

$${dE\over dk} = {E(k)- K(k)\over k}.\eqno(3.5)$$

It is known \cite{wit},  that $K(k),  K(k^\prime)$ satisfy  the differential equation

$${d\over dk } \left( k (k^\prime)^2 {du\over dk} \right)= k u \eqno(3.6)$$
and $E(k),\   E(k^\prime)- K(k^\prime)$ are, in turn,  solutions of  the differential equation

$$(k^\prime)^2 {d\over dk } \left( k  {du\over dk}\right) + ku =0.\eqno(3.7)$$  
The derivative of $K(k)$ can be calculated by the formula

$${dK\over dk} =   {E(k)-   (k^\prime)^2 K(k)\over k (k^\prime)^2} .\eqno(3.8)$$

Let $k_r$ be an elliptic modulus such that $x(k_r) = \sqrt r$ (see (3.1)).   In the sequel we will use such values for small $r$ and the corresponding elliptic integral singular values $K(k_r)$ (see \cite{bor}, \cite{borzuc}), namely

$$k_1= {1\over \sqrt 2}, \  k_2= \sqrt 2- 1,\  k_3= {1\over 4} \sqrt 2 (\sqrt 3-1),\ k_4= 3- 2\sqrt 2,\eqno(3.9)$$

$$K(k_1) = {\Gamma^2(1/4)\over 4\sqrt \pi},\quad     K(k_2) = {(\sqrt 2+1)^{1/2} \Gamma(1/8) \Gamma(3/8)\over 2^{13/4} \sqrt \pi},\eqno(3.10)$$

$$K(k_3) = {3^{1/4} \Gamma^3(1/3)\over 2^{7/3} \pi},\quad     K(k_4) = {(\sqrt 2+1) \Gamma^2 (1/4) \over 2^{7/2} \sqrt \pi}, \eqno(3.11)$$
where $\Gamma (z)$ is Euler's gamma-function \cite{erd}, Vol. I. According to \cite{borzuc} the so-called elliptic alpha function for the integral singular values 
$$\alpha(r)=  {E (k^\prime_r)\over K(k_r)} - {\pi\over 4 [K(k_r)]^2}=  {\pi\over 4 [K(k_r)]^2} + \sqrt r \left[ 1- 
{E (k_r)\over K(k_r)} \right] \eqno(3.12)$$
is calculated, in particular,  for small values and we have

$$\alpha(1)= {1\over 2},\   \alpha(2)=  \sqrt 2-1, \ \alpha(3)=  {1\over 2} (\sqrt 3-1),\  \alpha(4)=  2 (\sqrt 2-1)^2.\eqno(3.13)$$
Meanwhile,  appealing to relations (2.4.3.1), (2.4.3.3) in \cite{prud}  and the inverse Mellin transform \cite{tit}, we derive the following integral representations, related to the hyperbolic functions which will be useful in the sequel

$${1\over \sinh^2 ( cx)} = {2\over \pi i} \int_{\gamma-i\infty}^{\gamma+i\infty}  \Gamma(s) \zeta(s-1) (2cx)^{-s} ds, \ c >0, \gamma > 2, \eqno(3.14)$$

$${1\over \cosh^2( cx)} = {2\over \pi i} \int_{\gamma-i\infty}^{\gamma+i\infty}  (1- 2^{2-s}) \Gamma(s) \zeta(s-1 ) (2cx)^{-s} ds,\ c >0, \gamma > 0, \eqno(3.15)$$
where $\zeta(s)$ is the Riemann zeta-function \cite{erd}, Vol. I, which satisfies the  familiar functional equation
$$\zeta(s)= 2^s \pi^{s-1} \sin\left({\pi s\over 2}\right)\Gamma(1-s)\zeta(1-s).\eqno(3.16)$$

We begin, recalling the Rayleigh formula \eqref{eq:Ray} and recently obtained formula \eqref{eq:numT} in order to give a rigorous proof of the following functional equations for $S_2(\tau), T_2(\tau)$ on the imaginary positive half- axis and positive half-lines ${\rm Re}  \tau= \pm 1/2$.

{\bf Theorem 1}. {\it Let $x \in \mathbb{R}_+$. Then}

$$S_2(ix) +  S_2\left(i x^{-1}\right) = 2\pi,\eqno(3.17)$$

$$S_2\left( {\pm 1+ix\over 2} \right) +  S_2\left( {\pm 1+ix^{-1} \over 2} \right)$$

$$ = S_2\left( {\pm 1+ix\over 2} \right) +  S_2\left( {\mp 1+ix^{-1} \over 2} \right) = 2\pi,\eqno(3.18)$$

$$T_2(ix)= T_2\left(ix^{-1}\right), \eqno(3.19)$$

$$T_2\left( { \pm 1+ix\over 2} \right) -   T_2\left( { \pm 1+ix^{-1}\over 2} \right)$$

$$=  4 \left( S_2\left( { \pm 1+ix\over 2} \right) -  S_2 \left( ix \right) \right)  +  {2\pi^2 \over 3} \left(x-  {1\over x} \right). \eqno(3.20)$$

\begin{proof}  Indeed,  employing integral representation (3.14), we substitute it into \eqref{eq:Ray}  to write for $\tau= ix$

$$S_2(ix)=  \pi^2 x \left[ {1\over 3} -  {4\over \pi i} \sum_{m=1}^\infty  \int_{\gamma-i\infty}^{\gamma +i\infty}  \Gamma(s) \zeta(s-1) (2\pi m x)^{-s} ds\right].\eqno(3.21)$$
Since $\gamma > 2$ and the zeta- function is bounded on the vertical line $(\gamma-i\infty,\ \gamma + i\infty)$, i.e. $|\zeta(s-1)| \le \zeta (\gamma -1)$,  the interchange of the order of summation and integration is allowed for each $x >0$ via the absolute and uniform convergence by virtue of the estimate 

$$\sum_{m=1}^\infty  \int_{\gamma -i\infty}^{\gamma+i\infty}  \left| \Gamma(s) \zeta(s-1) (2\pi m x)^{-s} ds\right|$$

$$ \le  (2\pi x)^{-\gamma} \zeta(\gamma-1) \sum_{m=1}^\infty  {1\over m^\gamma } \int_{\gamma -i\infty}^{\gamma+i\infty}  \left| \Gamma(s) ds\right|  $$

$$=   (2\pi x)^{-\gamma} \zeta(\gamma-1) \zeta(\gamma)  \int_{\gamma-i\infty}^{\gamma+i\infty}  \left| \Gamma(s) ds\right|  < \infty,$$
where the convergence of the latter integral can be easily verified, appealing to the Stirling asymptotic formula for gamma-function when $|{\rm Im} s | \to \infty$ (see \cite{erd}, Vol. I).   Hence with the definition of the Riemann zeta-function in terms of the series, equality  (3.21) becomes

$$ S_2(ix)=  \pi^2 x \left[ {1\over 3} -  {4\over \pi i}  \int_{\gamma -i\infty}^{\gamma+i\infty}  \Gamma(s) \zeta(s) \zeta(s-1) (2\pi  x)^{-s} ds\right].\eqno(3.22)$$

On the other hand, the product of zeta-functions $\zeta(s) \zeta(s-1)$ can be represented by the Ramanujan identity \cite{yak}

$$\zeta(s) \zeta(s-1) = \sum_{m=1}^\infty {\sigma(m)\over m^s},\quad  \gamma > 2,\eqno(3.23)$$
where $\sigma(m)$ is the arithmetic function \cite{apo}, denoting the sum of divisors of $m$.  Hence, substituting in (3.22) and inverting the order of integration and summation owing to the same motivation, we obtain 

$$ S_2(ix)=  \pi^2 x \left[ {1\over 3} -  {4\over \pi i}  \sum_{m=1}^\infty \sigma(m) \int_{\gamma -i\infty}^{\gamma+i\infty}  \Gamma(s)  (2\pi m x)^{-s} ds\right]$$

$$=   \pi^2 x \left[ {1\over 3} -  8 \sum_{m=1}^\infty \sigma(m)   e^{-2\pi m x} \right],\eqno(3.24)$$
where the inverse Mellin transform of the gamma-function \cite{tit} is used

$$ e^{-x} =  {1\over 2\pi i} \int_{\gamma -i\infty}^{\gamma+i\infty}  \Gamma(s)  x^{-s} ds,\ x >0.$$

In the meantime,  the Nasim identity \cite{nasim} says that 

$$\sum_{m=1}^\infty \sigma(m)  e^{-2\pi m x}  + x^{-2} \sum_{m=1}^\infty \sigma(m)   e^{-2\pi m / x} =
{1\over 24} \left( 1+ {1\over x^2} \right) - {1\over 4\pi x},\  x > 0.\eqno(3.25)$$
Consequently, from (3.22) we find 

$$S_2(ix) +  S_2\left(i x^{-1}\right) =   {\pi^2\over 3} \left(  x + {1\over x} \right) - 8\pi^2 \left[ x \sum_{m=1}^\infty \sigma(m)   e^{-2\pi m x}\right.$$

$$\left.   +  {1 \over  x}   \sum_{m=1}^\infty \sigma(m)   e^{-2\pi m/ x} \right] = 2\pi,$$
proving equation (3.17).  In order to prove equations (3.18), we invoke representation  (3.15), motivating all passages analogously to the previous case. Moreover, as we will see it is sufficient to prove (3.18) for positive real parts. So, we have  (see \eqref{eq:Ray}) 

$$S_2\left( { 1+ix\over 2} \right) =  {\pi^2 x\over 2}  \left[ {1\over 3}  -  2 \sum_{m=1}^\infty  {1\over \sinh^2( \pi m x)} +
2 \sum_{m=1}^\infty  {1\over \cosh^2( \pi (m - 1/2) x)} \right]$$

$$=  {1\over 2}\ S_2(ix) -   2 \pi i  x   \int_{\gamma-i\infty}^{\gamma+i\infty}  (1- 2^{2-s}) \Gamma(s) \zeta(s-1 ) (\pi x)^{-s}
\sum_{m=1}^\infty {1\over (2m-1)^s } ds. $$
But the latter series is easily calculated for $\gamma > 1$ via the definition of the Riemann zeta-function and we obtain

$$\sum_{m=1}^\infty {1\over (2m-1)^s } =  (1- 2^{-s} ) \zeta (s).$$
Hence, recalling the Ramanujan identity (3.23),  changing the order of integration and summation and calculating the inverse Mellin transform of the gamma-function of different arguments, we deduce 

$$S_2\left( { 1+ix\over 2} \right) =   {1\over 2}\  S_2(ix) -   2 \pi i  x   \int_{\gamma-i\infty}^{\gamma+i\infty}  (1- 2^{2-s})  (1- 2^{-s} ) $$

$$\times \Gamma(s)  \zeta (s) \zeta(s-1 ) (\pi x)^{-s} ds=  {1\over 2}\  S_2(ix)$$

$$+ 4\pi^2 x \sum_{m=1}^\infty \sigma(m) \left[ e^{- \pi m x} -  5 e^{- 2 \pi m x} + 4 e^{- 4 \pi m x} \right] .$$
Meanwhile, from (3.24) we find

$$x \sum_{m=1}^\infty \sigma(m) e^{- 2\pi m x} = { x\over 24}  - {1\over 8\pi^2}\  S_2(ix).\eqno(3.26)$$
Therefore, it yields

$$S_2\left( { 1+ix\over 2} \right) =    S_2(ix)  - {\pi^2 x\over 6} + 4\pi^2 x \sum_{m=1}^\infty \sigma(m) \left[ e^{- \pi m x} -  2 e^{- 2 \pi m x}\right ] $$

$$-  4\pi^2 x \sum_{m=1}^\infty \sigma(m) \left[ 2 e^{- 2\pi m x}- 4 e^{- 4 \pi m x} \right].\eqno(3.27)$$
In the meantime, appealing to another Nasim's  formula (see \cite{nasim}, formula (5.1) with $a= x/2,\ b=x $)

$$ \sum_{m=1}^\infty {\sigma(m) \over m} \left[ e^{- \pi m x} -   e^{- 2 \pi m x} \right]  = \sum_{m=1}^\infty {\sigma(m) \over m} \left[ e^{- 4\pi m/ x} -   e^{- 2 \pi m/ x} \right]$$

$$  +  {\pi\over 12} \left( {1\over x} + {x\over 2} \right) - {1\over 2}\log 2,$$
we differentiate it with respect to $x$, which is permitted via the absolute and uniform convergence and  multiply by $x$ the obtained equality. Thus we get

$$ x \sum_{m=1}^\infty \sigma(m) \left[ 2 e^{- 2\pi m x} -   e^{- \pi m x} \right]  = {2\over x} \sum_{m=1}^\infty \sigma(m) \left[ 2 e^{- 4\pi m/ x} -   e^{- 2 \pi m/ x} \right] $$

$$+ {x\over 24}   -   {1\over 12 x}.\eqno(3.28)$$
Substituting in (3.27),  we derive

$$S_2\left( { 1+ix\over 2} \right) =    S_2(ix)    + { 8\pi^2\over  x} \sum_{m=1}^\infty \sigma(m) \left[  e^{- 2 \pi m/ x}- 2 e^{- 4\pi m/ x}    \right]  $$

$$-  8\pi^2 x \sum_{m=1}^\infty \sigma(m) \left[  e^{- 2\pi m x}- 2 e^{- 4 \pi m x} \right].\eqno(3.29)$$
Now, changing in  (3.29)  $x$ by $1/x$ and adding these two equalities with the use of (3.17),  we obtain  (3.18).

Let us prove (3.19).   To do this,   we let $\tau= ix,\ x >0$ in \eqref{eq:numT} and write it  in the form

$$T_2(ix)= S_2(ix) + 4\pi^3 x^2 \sum\limits_{n=1}^\infty n \frac{\cosh( n\pi x)}{\sinh^3 (n\pi x)}.\eqno(3.30) $$

However, the series in (3.30) can be obtained by termwise differentiation with respect to $x$ of the series $\sum \hbox{cosech}^2(\pi n x)$ for $x \ge x_0 >0$ due to the absolute and uniform convergence.   Hence

$$T_2(ix)= S_2(ix) -  2\pi^2 x^2 {d\over dx} \sum\limits_{n=1}^\infty  \frac{1}{\sinh^2 (n\pi x)}.\eqno(3.31) $$
But from \eqref{eq:Ray}, (3.23), (3.24)   and termwise differentiation of the series with arithmetic function $\sigma(m)$ in (3.24) by the same reasons, we obtain

$$T_2(ix)= S_2(ix) + 16\pi^3 x^2  \sum\limits_{n=1}^\infty  m\  \sigma(m) e^{-2\pi m x}.\eqno(3.32) $$
Meanwhile, differentiating the Nasim identity (3.25) with respect to $x$ and then multiplying both sides of the obtained equality by $- x^2/ (2\pi) $, we find

$$x^2 \sum_{m=1}^\infty m\  \sigma(m)  e^{-2\pi m x} =  x^{-2} \sum_{m=1}^\infty m\  \sigma(m)   e^{-2\pi m / x} +
{1\over 24 \pi x}$$

$$  - {1\over 8\pi^2}- {1\over \pi x} \sum_{m=1}^\infty  \sigma(m)   e^{-2\pi m / x}.\eqno(3.33)$$
Substituting the left-hand side of the latter equality in (3.32) and appealing to (3.26), we write 

$$T_2(ix)-  S_2(ix) =  T_2\left(ix^{-1}\right) + S_2\left(ix^{-1}\right) - 2\pi.$$
Therefore, equality (3.17) leads us to (3.19).     In order to establish (3.20), we employ again \eqref{eq:numT} to get as in (3.31), (3.32) 

$$T_2\left( { \pm 1+ix\over 2} \right) =  S_2 \left( { \pm 1+ix\over 2} \right)  + (\pi  x)^2 {d\over dx} \left[  \sum_{m=1}^\infty  {1 \over \cosh^2( \pi (m - 1/2) x)}  \right.$$

$$\left.  -   \sum_{m=1}^\infty  {1\over \sinh^2( \pi m x)} \right]= S_2 \left( { \pm 1+ix\over 2} \right)  + \pi^3  x^2  \left[  {1\over \pi} {d\over dx} \sum_{m=1}^\infty  {1 \over \cosh^2( \pi (m - 1/2) x)}  \right.$$

$$\left.  + 8   \sum\limits_{n=1}^\infty  m\  \sigma(m) e^{-2\pi m x} \right].$$
In the meantime, recalling (3.15), (3.23) and termwise differentiation, we deduce

$${d\over dx} \sum_{m=1}^\infty  {1 \over \cosh^2( \pi (m - 1/2) x)} = 4 \pi  \sum\limits_{n=1}^\infty m\  \sigma(m) \left[ 10 e^{-2 \pi m x}\right.$$

$$\left.  - 16 e^{-4 \pi m x}  - e^{- \pi m x}\right] . $$
Thus

$$T_2\left( { \pm 1+ix\over 2} \right) =  S_2 \left( { \pm 1+ix\over 2} \right)  +  4 \pi^3  x^2 \sum\limits_{n=1}^\infty m\  \sigma(m) \left[ 12 e^{-2 \pi m x}\right.$$

$$\left.  - 16 e^{-4 \pi m x}  - e^{- \pi m x}\right].$$
Hence,

$${1\over 4\pi^3} \left(  T_2\left( { \pm 1+ix\over 2} \right) - S_2 \left( { \pm 1+ix\over 2} \right) \right) =
8 x^2 \sum\limits_{n=1}^\infty m\  \sigma(m)  e^{-2 \pi m x} $$

$$-  16 x^2 \sum\limits_{n=1}^\infty m\  \sigma(m) e^{-4 \pi m x}
- x^2 \sum\limits_{n=1}^\infty m\  \sigma(m) \left[ e^{- \pi m x}- 4  e^{-2 \pi m x}\right] .\eqno(3.34)$$
Returning to (3.28) and making the termwise differentiation and simple manipulations, we  derive

$$  x^2 \sum_{m=1}^\infty m\ \sigma(m) \left[  e^{- \pi m x} - 4 e^{- 2\pi m x}  \right]  = -  {4\over \pi x} \sum_{m=1}^\infty \sigma(m) \left[ 2 e^{- 4\pi m/ x} -   e^{- 2 \pi m/ x} \right] $$

$$    +  {1\over 6 \pi x} +   {4\over x^2} \sum_{m=1}^\infty m\ \sigma(m) \left[ 4 e^{- 4\pi m/ x} -   e^{- 2 \pi m/ x} \right].$$
Therefore with the use of (3.32), equality (3.34) becomes

$$ {1\over 4\pi^3} \left(  T_2\left( { \pm 1+ix\over 2} \right) - S_2 \left( { \pm 1+ix\over 2} \right) \right) =
{1\over 2\pi^3} \left(  T_2\left( ix \right) - S_2 \left( ix \right) \right)  $$

$$-  16 x^2 \sum\limits_{n=1}^\infty m\  \sigma(m) e^{-4 \pi m x} + { 4\over \pi x} \sum_{m=1}^\infty \sigma(m) \left[  e^{- 2 \pi m/ x}- 2 e^{- 4\pi m/ x}  \right] 
$$

$$- {1\over 6\pi x}  -  {16\over  x^2} \sum\limits_{n=1}^\infty m\  \sigma(m) e^{-4 \pi m/ x} +  {1\over 4\pi^3} \left(  T_2\left( ix^{-1} \right) - S_2 \left( ix^{-1} \right) \right).$$  
Hence,

$$ {1\over 4\pi^3} \left(  T_2\left( { \pm 1+ix\over 2} \right) - S_2 \left( { \pm 1+ix\over 2} \right) \right) -
{1\over 2\pi^3} \left(  T_2\left( ix \right) - S_2 \left( ix \right) \right)  $$

$$-  {1\over 4\pi^3} \left(  T_2\left( ix^{-1} \right) - S_2 \left( ix^{-1} \right) \right) +  {1\over 6\pi x} - { 4\over \pi x} \sum_{m=1}^\infty \sigma(m) \left[  e^{- 2 \pi m/ x}- 2 e^{- 4\pi m/ x}  \right]$$

$$=  {1\over 4\pi^3} \left(  T_2\left( { \pm 1+ix^{-1}\over 2} \right) - S_2 \left( { \pm 1+ix^{-1} \over 2} \right) \right) -
{1\over 2\pi^3} \left(  T_2\left( ix^{-1} \right) - S_2 \left( ix^{-1} \right) \right)  $$

$$-  {1\over 4\pi^3} \left(  T_2\left( ix \right) - S_2 \left( ix\right) \right) +  {x\over 6\pi} - { 4x\over \pi } \sum_{m=1}^\infty \sigma(m) \left[  e^{- 2 \pi m x}- 2 e^{- 4\pi m x}  \right].$$
Meanwhile,  appealing to (3.29),  we find 

$$ { 4\over \pi x} \sum_{m=1}^\infty \sigma(m) \left[  e^{- 2 \pi m/ x}- 2 e^{- 4\pi m/ x}    \right]  -  {4x\over \pi} \sum_{m=1}^\infty \sigma(m) \left[  e^{- 2\pi m x}- 2 e^{- 4 \pi m x} \right] $$

$$= {1\over 2\pi^3} \left( S_2\left( { 1+ix\over 2} \right) -    S_2(ix)\right).$$
Thus, accounting (3.19), 

$$   T_2\left( { \pm 1+ix\over 2} \right) -   T_2\left( { \pm 1+ix^{-1}\over 2} \right) +3 S_2 \left( ix \right) -  S_2 \left( ix^{-1} \right) $$

$$=  3\  S_2\left( { \pm 1+ix\over 2} \right) - S_2 \left( { \pm 1+ix^{-1} \over 2} \right)  +  {2\pi^2 \over 3} \left(x-  {1\over x} \right).$$
Finally, equalities (3.17), (3.18) drive us at (3.20),   completing the proof of Theorem 1. 
\end{proof} 

The explicit expressions of $S_2(\tau)$ on the imaginary axis and the lines ${\rm Re} \tau= \pm 1/2$ are given by 

{\bf Theorem 2}.   {\it Let  $x \in \mathbb{R} \backslash\{0\}$. Then the following formulae hold

$$S_2(ix)=  {4\over 3}\  \hbox{sign}( x)  \  K(k^\prime) \left[ 3E(k) + (k^2- 2) K(k) \right],\eqno(3.35)$$

$$S_2\left( {\pm 1+ix\over 2} \right) =  2\  \hbox{sign}( x) \   K(k^\prime) \left[ 2\  E(k)+  { 4k^2 -5\over 3} \  K(k)\right], \eqno(3.36)$$
where  

$$|x |= { K (k^\prime)\over K(k)},\  k \in (0, 1)$$
and $k^\prime$ is defined by $(3.1)$.} 

\begin{proof}  Let us first consider a positive $x$ being  defined by (3.1).   Fortunately,  the series in \eqref{eq:Ray} for $\tau= ix$ is calculated in \cite{prud}, relation (5.3.4.5),  and we have

$$\sum_{m=1}^\infty {1\over \sinh^2(\pi m x )}=  {1\over 6}+ {2(2-k^2) \over 3\pi^2}\ K^2(k) - {2\over \pi^2} K(k)E(k).\eqno(3.37)$$ 
Therefore, 

$$S_2(ix)= {4\over 3} \ (k^2- 2) \ x K^2(k)  + 4x  K(k)E(k) $$

$$=  {4\over 3}\  K(k^\prime) \left[ 3E(k) + (k^2- 2) K(k) \right].$$
Hence it proves (3.35) for positive $x$, and for negative $x$ it can be easily extended via \eqref{eq:Ray}.  In order to prove (3.36), we employ relation (5.3.6.6) in \cite{prud}

$$\sum_{m=1}^\infty {1\over \cosh^2(\pi x (m-1/2) )}=  {2\over \pi^2} K(k)E(k) -   {2(1-k^2) \over \pi^2}\ K^2(k).\eqno(3.38)$$ 
Then for positive $x$ we find from \eqref{eq:Ray}, (3.37), (3.38) 

$$S_2\left( {\pm 1+ix\over 2} \right) =  \pi^2 x  \left[\frac 16+
\sum _{m=1}^{\infty } \frac{1}{\sin^2(\pi  m (\pm 1+ ix))} \right.$$

$$\left. + \sum _{m=1}^{\infty } \frac{1}{\sin^2(\pi  (2m-1) (\pm 1+ ix)/2)} \right] = \pi^2 x  \left[\frac 16- 
\sum _{m=1}^{\infty } \frac{1}{\sinh^2(\pi  m x)} \right.$$

$$\left. +  \sum _{m=1}^{\infty } \frac{1}{\cosh^2 (\pi x (m-1/2))} \right] =   4 K(k^\prime)E(k) - {2(2-k^2) \over 3}\ K(k^\prime) K(k)$$

$$-  2(1-k^2) \  K(k^\prime) K(k) =  2 K(k^\prime) \left[ 2 E(k) + \frac{4k^2 -5}{3}\  K(k)\right].$$
Hence spreading  the latter equalities for negative $x$, we get (3.36).
\end{proof}  

{\bf Corollary 1}. {\it Formula $(3.20)$ can be written in the form }

$$T_2\left( { \pm 1+ix\over 2} \right) -   T_2\left( { \pm 1+ix^{-1}\over 2} \right)=   {4\over 3}  \left( 4k^2 -2\right) K(k)  +  {2\pi^2 \over 3} \left(x-  {1\over x} \right).$$

As we could see above,  the only value $S_2(i)=\pi$ was known explicitly. Now we are able to calculate more  interesting particular values of (3.35), (3.36).  Indeed, we have  

{\bf Corollary 2}.   {\it The following values take place} 

$$ S_2(\pm i)= S_2\left( {1\pm i\over 2} \right) = \pm \  \pi,\eqno(3.39)$$

$$S_2(\pm i\sqrt 2) = \pm \left[ \pi +  {\Gamma^2(1/8) \Gamma^2(3/8)\over 48  \pi\sqrt 2}\right],\eqno(3.40)$$

$$S_2\left( {1\pm i\sqrt 2 \over 2} \right) = \pm  \left[ \pi + \frac { (2\sqrt 2- 3)\   \Gamma^2(1/8)\Gamma^2(3/8) } { 96 \ \pi} \right],\eqno(3.41)$$ 

$$S_2(\pm i\sqrt 3 )= \pm \left[  \pi + {\sqrt 3\  \Gamma^6( 1/3) \over 16  \pi^2\   2^{2/3}} \right],\eqno(3.42)$$

$$S_2\left( {1\pm i\sqrt 3 \over 2} \right)=  \pm \  \pi,\eqno(3.43)$$

$$S_2\left(\pm 2i\right) = \pm \left[ \pi +  {\Gamma^4(1/4) \over 16\ \pi} \right], \eqno(3.44)$$

$$S_2\left( {{1\over 2} \pm i} \right) = \pm  \left[ \pi + \frac { (3- 2\sqrt 2)\   \Gamma^4 (1/4) } { 32 \ \pi} \right].\eqno(3.45)$$ 

\begin{proof}  As we observe from (3.35), (3.36),  it is sufficient to establish the above constants for a positive imaginary part of the corresponding $\tau$.  To do this we employ particular cases (3.9) of the modulus $k_r$  and the corresponding singular values (3.10), (3.11) $K(k_r),\ r= 1,2,3,4$.    In fact, letting $x=1,\ \sqrt 2,\ \sqrt 3,\  2$ and taking in mind (3.12), (3.13), we derive, respectively, 

$$S_2(i)= {4\over 3} K(k_1) \left[ {3\over 2} K(k_1) + {3\pi\over 4 K(k_1)} - {3\over 2} K(k_1) \right] = \pi;$$

$$S_2\left( {1+ i\over 2} \right) =  2 K(k_1) \left[ K(k_1) + {\pi\over 2 K(k_1)} - K(k_1) \right] = \pi;$$

$$S_2(i\sqrt 2)= {4\sqrt 2 \over 3} \left[ {3\over \sqrt 2} K^2(k_2) + {3\pi\over 4\sqrt 2} + (1-2\sqrt 2) K^2(k_2) \right] $$
$$ = \pi + {4\over 3} \left( \sqrt 2 - 1 \right) K^2(k_2) = \pi +  {\Gamma^2(1/8) \Gamma^2(3/8)\over 48\  \pi \sqrt 2 } ;$$

$$S_2\left( {1+  i\sqrt 2 \over 2} \right) = \pi + {2(7\sqrt 2-10) \over 3} \  K^2(k_2) = \pi + \frac { (2\sqrt 2- 3)\   \Gamma^2(1/8)\Gamma^2(3/8) } { 96 \ \pi} ;$$

$$S_2( i\sqrt 3 )=  \pi +  K^2(k_3)  =   \pi + {\sqrt 3\  \Gamma^6( 1/3) \over 16  \pi^2\   2^{2/3}}; $$

$$S_2\left( {1+  i\sqrt 3 \over 2} \right)=   \pi + 2K^2(k_3) (\sqrt 3 +1) - {2\over \sqrt 3} (3+\sqrt 3 ) K^2(k_3) = \pi;$$

$$S_2( 2i)= \pi +  8 K^2(k_4) (3-2\sqrt 2) = \pi +  {\Gamma^4(1/4) \over 16\ \pi};$$

$$S_2\left( {{1\over 2} + i} \right) =  \pi + 4(17-12\sqrt 2) K^2(k_4) =  \pi +  \frac { (3- 2\sqrt 2)\   \Gamma^4 (1/4) } { 32 \ \pi}.$$
\end{proof}  

A  more technically difficult task is to find explicit expressions for $T_2(\tau)$ on the same lines in the complex plane.  To achieve our goal we will involve the method of termwise differentiation of the series in \eqref{eq:Ray} with respect to the elliptic modulus (for $\tau= ix(k)$ or $\tau = (\pm 1 + i x(k))/2 )$.   Indeed,  as we mentioned above,  $x(k)$ by formula (3.1) is continuously differentiable and when $k \in [a_0,  b_0],\  0< a_0 <  b_0 <1$,  the corresponding series \eqref{eq:Ray} is absolutely and uniformly convergent. Moreover, it is not difficult to show that  the series of the derivatives with respect to $k$ converges absolutely and uniformly on the segment $[a_0,  b_0]$. Thus the known theorem from calculus says that the termwise differentiation of the series is allowed.    This leads us to 

{\bf Theorem 3}.   {\it Under conditions of Theorem $2$  the following formulae hold valid}

$$T_2(ix) =  \frac 43 \   \hbox{sign} ( x) \  K(k') \left[   \left[ 1-  {2 \over \pi } K(k') E(k) \right]  
\left[ 3E(k) + (k^2- 2) K(k) \right]\right.$$

$$\left. - {2\over \pi} K(k^\prime) K(k) \left[ (1-k^2) \left[K(k)- E(k)\right] - E(k)\right]\right],\eqno(3.46)$$

$$T_2\left( { \pm 1+ix\over 2} \right) =   {2\over 3}\  \hbox{sign} ( x)   K(k^\prime) \left[  \left( 6 E(k) +  K(k) (4k^2-5) \right)\right.$$

$$\left. \times  \left[ 1 - {2\over \pi} K(k^\prime) \left(  E(k) +  K(k) (k^2-1) \right)\right]\right. $$

$$\left.  -  {2\over \pi}   K(k^\prime) K(k) \left[  (1-2k^2)  E(k)+   (4k^2-1) (1-k^2) K(k) \right] \right].\eqno(3.47)$$

\begin{proof}  Indeed, concerning the proof of formula (3.37), we let $\tau= ix,\ x >0$ in \eqref{eq:numT} and write it  in the form

$$T_2(ix)= S_2(ix) + 4\pi^3 x^2 \sum\limits_{n=1}^\infty n \frac{\cosh( n\pi x)}{\sinh^3 (n\pi x)},\eqno(3.48) $$
where $x$ is a function of $k$ by (3.1) and since the termwise differentiation is permitted, we obtain

$$\sum\limits_{n=1}^\infty n \frac{\cosh( n\pi x)}{\sinh^3 (n\pi x)}  = - { 1\over 2\pi x^\prime(k)} {d\over dk} \sum\limits_{n=1}^\infty  \frac{1}{\sinh^2 (n\pi x)}.\eqno(3.49)$$

Meanwhile,  the derivative $x^\prime(k)$ can be calculated explicitly, employing twice (3.8). Hence  we find

$$x^\prime (k) =  -  \frac{K(k^\prime)} {K^2(k)} {dK(k)\over dk} -  \frac{k } {(1-k^2) K(k)} \left[ \frac{E (k^\prime)} {k^2} -
K(k^\prime)\right]$$

$$=     \frac{1} {k(1-k^2) K(k)} \left[ K(k^\prime) \left[ 1 -  \frac{E(k)} {K(k)} \right] -   E (k^\prime)\right] $$ 
and the Legendre identity (3.3) leads  us  to the final result

$$x^\prime (k) =  - \frac{\pi}{2 k(1-k^2) K^2(k)}.\eqno(3.50)$$
Therefore, recalling (3.5), (3.8), (3.37), we deduce from (3.49)

$$4\pi^3 x^2 \sum\limits_{n=1}^\infty n \frac{\cosh( n\pi x)}{\sinh^3 (n\pi x)} =  {8\over \pi}  k(1-k^2) K^2(k^\prime) {d\over dk} \left[  K(k) \left( {2-k^2\over 3} K(k) - E(k) \right) \right]$$

$$= {8\over 3\pi} K^2(k^\prime) \left( E(k) - (1-k^2) K(k) \right) \left(  (2-k^2) K(k) - 3 E(k) \right)$$

$$+  {8\over 3\pi}   K^2(k^\prime) K(k) \left[  E(k) (2k^2-1) + K(k) (k^2-1)^2  \right]$$

$$=   {8\over 3\pi} K^2(k^\prime) \left[ 2 E(k)K(k)  (2- k^2) + K^2(k) (k^2-1) - 3 E^2(k) \right].$$
Hence, appealing to (3.35) and combining with (3.48), we arrive at (3.46) being valued for positive $x$. Then we extend it on negative numbers as in Theorem 2. 

In order to establish identity (3.47), we write \eqref{eq:numT} for $\tau= (\pm 1+ ix)/2,\ x >0$ in the same manner as in the proof of identity (3.20).    Nevertheless, we will employ explicit expressions (3.37) and (3.38) and make the termwise differentiation with respect to the elliptic modulus. Hence,  taking in mind (3.36), (3.50),  we obtain 

$$T_2\left( { \pm 1+ix\over 2} \right) =  S_2 \left( { \pm 1+ix\over 2} \right)  + {(\pi  x)^2 \over x^\prime (k)} {d\over dk} \left[  \sum_{m=1}^\infty  {1 \over \cosh^2( \pi (m - 1/2) x)}  \right.$$

$$\left.  -   \sum_{m=1}^\infty  {1\over \sinh^2( \pi m x)} \right] = {2\over 3} K(k^\prime)  \left( 6 E(k) +  K(k) (4k^2-5) \right)$$

$$  -  {4\over 3 \pi}  k(k^\prime)^2  K^2(k^\prime)  {d\over dk} \left[ K(k) \left( 6 E(k) +  K(k) (4k^2-5) \right) \right].$$

Fulfilling the differentiation with the aid of  (3.5), (3.6), (3.8), we find 

$$T_2\left( { \pm 1+ix\over 2} \right) =   {2\over 3} K(k^\prime) \left[  \left( 6 E(k) +  K(k) (4k^2-5) \right)\right.$$

$$\left. \times  \left[ 1 - {2\over \pi} K(k^\prime) \left(  E(k) +  K(k) (k^2-1) \right)\right]\right. $$

$$\left.  -  {2\over \pi}   K(k^\prime) K(k) \left[  (1-2k^2)  E(k)+   (4k^2-1) (1-k^2) K(k) \right] \right], $$
arriving at (3.47) after the same extension on negative numbers $x$ as in Theorem 3.
\end{proof}

As a corollary we calculate particular values of $T_2$ on the mentioned vertical lines, recalling $k_r$ in  (3.9) and $K(k_r)$ in  (3.10), (3.11), letting $\ r= 1,2,3,4$.   In particular, the value $x=1$, corresponding $k_1=k_1'=\frac{1}{\sqrt{2}}$, gives  the interesting and important constant numerical value of which coincides with the numerical value of $T_2(i)= 4.078451$ computed with \eqref{eq:numT}  
$$T_2(i)=  {\pi\over 2} + {\Gamma^8(1/4) \over 384\  \pi^3}.\eqno(3.51)$$
We note that this  numerical result $T_2(i)= 4.078451$ coincides with the numerical value obtained by other approaches \cite{McPhedran}, \cite{Greengard}. 

{\bf Corollary 3}. {\it Certain explicit constants related to $T_2(\tau)$ are the following values}

$$T_2(\pm i) = \pm \left[ {\pi\over 2} + {\Gamma^8(1/4) \over 384\  \pi^3}\right] ; $$ 

$$T_2\left( {  1 \pm i\over 2} \right) = \pm  \left[ {\pi\over 2} - {\Gamma^8(1/4) \over 384\  \pi^3}\right] ; $$ 

$$T_2(\pm i\sqrt 2) = \pm \left[ {\pi\over 2} + {\Gamma^4(1/8)\Gamma^4(3/8) \over 1024\  \pi^3} \right]; $$ 

$$T_2\left( {  1\pm i\sqrt 2 \over 2} \right) = \pm  \left[ {\pi\over 2} - {\Gamma^4(1/8)\Gamma^4(3/8) (\sqrt 2 -1) \over 1024\  \pi^3} \right]; $$  

$$T_2\left( {  \pm i\sqrt 3} \right) = \pm  \left[  {\pi\over 2} -  {2^{2/3} \Gamma^{12} (1/3) (9+ 4\sqrt 3) \over 512\  \pi^5} \right]; $$

$$T_2\left({1\pm i\sqrt 3\over 2}\right) = \pm  {\pi\over 2}; $$

$$T_2\left( {  \pm 2i} \right) = \pm  \left[  {\pi\over 2} +   { \Gamma^{8} (1/4) \over 192\  \pi^3} \right]; $$

$$T_2\left({1\over 2} \pm i \right) = \pm   \left[  {\pi\over 2} +  {\Gamma^{8} (1/4) (5-3\sqrt 2) \over 768\  \pi^3} \right]; $$

\section{Random lattice sums}
\label{sec:5}
Consider a lattice with the fixed periods $\omega_1$, $\omega_2$ and the corresponding fundamental parallelogram
$$
\mathcal{G}_{(0,0)}:=\left\{t_1\omega_1 + t_2\omega_2 \in \mathbb{C}: 0<t_1,t_2<1\right\}
$$
The Eisenstein function of second order \cite{Weil} is defined by the series
\begin{equation}
\label{eq:ran1}
E_2(z) =\sideset{}{^e}\sum_{(m,n) \in \mathbb Z^2}
\frac{1}{(z-m \omega_1-n \omega_2)^2}.
\end{equation}
It is related to the $\wp$-Weierstrass function by formula \cite{Weil}
\begin{equation}
\label{eq:ran2}
E_2(z) =\wp(z) +S_2.
\end{equation} 
Following \eqref{eq:ran1} we introduce the function 
\begin{equation}
\label{eq:ran3}
G_2(z) =\sideset{}{^e}\sum_{(m,n) \in \mathbb Z^2}
\frac{\overline{z-m \omega_1-n \omega_2}}{(z-m \omega_1-n \omega_2)^3}.
\end{equation}
This function is related to the Natanzon function \cite{Natanzon}
\begin{equation}
\label{eq:ran4}
\wp'_1(z) =-2 \sum_{(m,n) \in \mathbb Z^2 \backslash \{(0,0)\}}
\left[
\frac{\overline{z-m \omega_1-n \omega_2}}{(z-m \omega_1-n \omega_2)^3} +
\frac{\overline{m \omega_1+n \omega_2}}{(m \omega_1+n \omega_2)^3}
\right]
\end{equation}
by formula
\begin{equation}
\label{eq:ran5}
G_2(z) = -\frac 12 \overline{z} \wp'(z)+\frac 12\wp'_1(z) + T_2.
\end{equation}
Filshtinsky \cite[Appendix 2]{Fil1992} found a relation between the Natanzon and Weierstrass functions which can be written in our case as  
\begin{equation}
\label{eq:ran6}
\pi \wp'_1(z) = \frac 13 \wp''(z)+[\zeta(z)-(S_2-\pi)z]\wp'(z)- 2(S_2-\pi)\wp(z) - 10 S_4,
\end{equation}
where $\zeta(z)$ is the $\zeta$-Weierstrass function and $S_4$ is defined by the absolutely convergent series 
\begin{equation}
\label{eq:ranS4}
S_4 =\sum_{(m,n) \in \mathbb Z^2 \backslash \{(0,0)\}}
\frac{1}{(m \omega_1+n \omega_2)^4}.
\end{equation}
Substitution of \eqref{eq:ran6} into \eqref{eq:ran5} yields
\begin{eqnarray}
\nonumber
G_2(z) =-\frac{1}2\overline{z}\wp'(z)+\frac{1}{6\pi} \wp''(z)+
\frac{1}2 \left[\frac{\zeta(z)}{\pi} -\left(\frac{S_2}{\pi}-1 \right) z \right]\wp'(z)
\\
-\left(\frac{S_2}{\pi}-1 \right) \wp(z) -\frac{5}{\pi} S_4+T_2.
\label{eq:e32}
\end{eqnarray}

Consider $N$ non-overlapping circular disks $D_k$ of radius $r$ with the centers $a_k\in\mathcal{G}_{(0,0)}$. These centers can be considered as random variables. 
Introduce the sums
\begin{equation}
\label{eq:ran8}
e_2 = \frac{1}{N^2} \sum_{k=1}^N \sum_{m=1}^N E_2(a_k-a_m),
\end{equation}
\begin{equation}
\label{eq:ran9}
g_2 = \frac{1}{N^2} \sum_{k=1}^N \sum_{m=1}^N F_2(a_k-a_m),
\end{equation}
where it is assumed that $E_2(0):=S_2$ and $F_2(0):=T_2$. Such a consideration implies that for $N=1$ $e_2$ becomes $S_2$ and $g_2$ becomes $T_2$. 

The sums \eqref{eq:ran8}-\eqref{eq:ran9} play the fundamental role in the theory of random 2D composites, since the effective conductivity tensor of the composite represented by $N$ discs per periodicity cell can be calculated by the asymptotic formula \cite{MR2001} (cf. \eqref{eq:CMA})
\begin{align}
\label{eq:CMAa}
\frac{\lambda_{xx}-i \lambda_{xy}}{\lambda} &= 1+2\rho f+2\rho^2 f^2 \frac{e_2}{\pi} + O((|\rho| f)^3),\\
\frac{\lambda_{yy}+i \lambda_{xy}}{\lambda} &= 1+2\rho f+2\rho^2 f^2 \left(2-\frac{e_2}{\pi} \right)+ O((|\rho| f)^3),
\end{align}
In the case of macroscopically isotropic composites, $\lambda_{xx}=\lambda_{yy}$ and $\lambda_{xy}=0$. This implies that $e_2$ must be equal to $\pi$. One can consider this assertion as a physical prove of the identity $e_2 = \pi$ for a macroscopically isotropic distribution of $a_k$.

Analogous formulae take place for the elastic constants. Let elastic fibers $D_k$  with the shear modulus $\mu_1$ and the Poisson ration $\nu_1$ are distributed in the matrix with the constants $\mu$ and $\nu$. Let $\kappa = 3-4\nu$ and $\kappa_1 = 3-4\nu_1$ be the corresponding Muskhelishvili constants for the plane strain. Consider the averaged constant $\mu_{e}= \frac{\langle \sigma_{xx}-\sigma_{yy} \rangle}{2\langle \epsilon_{xx} - \epsilon_{yy} \rangle}$ where $\sigma_{\alpha \beta}$ and $\epsilon_{\alpha \beta}$ denote the components of the stress and deformation tensors, respectively ($\alpha$ and $\beta$ can be $x$ and $y$). Here, $\langle \cdot \rangle$ denotes the average value (double integral over the periodicity cell). In particular, for macroscopically isotropic composites $\mu_{e}$ yields the effective shear modulus. The value $\mu_{e}$ can be estimated by asymptotic formula deduced in \cite{arx}  
\begin{equation}
\frac{\mu_{e}}{\mu (1+\kappa)}= \frac{1}{1+\kappa} + \frac{{\mu_1}-{\mu}}{\kappa{\mu_1}+{\mu} } f +\left(\frac{{\mu_1}-{\mu}}{\kappa{\mu_1}+{\mu} } \right)^2 \left(\kappa  -\frac{2\mbox{Re}\;g_2}{\pi} \right) f^2  +O(f^3).
\label{eq:eff}
\end{equation}
One can see that the value $g_2$ from \eqref{eq:ran9} occurs in the coefficient on $f^2$. 

Numerical simulations of  $e_2$ were performed in \cite{MNCz2012} for macroscopically isotropic composites generated by the RSA algorithm and by random walks. Using the RSA protocol we compute $100$ times $g_2$ for $r=0.003$ when $f$ is about $0.09$, $N$ is about $3250$. More precisely, $f$ and $N$ slightly change in each simulation of location in accordance with the RSA protocol \cite{MNCz2012}.  The mean value of $g_2- \frac{\pi}{2}$ holds $0.00457824 + 0.0121335 i$, the variance $0.0286453$. For $e_2- \pi$ we get the mean value $0.000723263 + 0.00575626 i$ and the variance $	0.0296968$.

\section{Conclusion}
\label{sec:6}
Explicit formulae of Section 3  deduced in this paper yield asymptotic analytical formulae for the effective tensors of 2D composites with circular inclusions. The obtained fundamental values $S_2$ and $T_2$ give a possibility to pass through the approximation $O(f^2)$ terms to get high order analytical formulae for the effective elastic constants of fibrous composites \cite{arx}.

\end{document}